\documentclass[12pt]{article}
\usepackage[intlimits]{amsmath}
\usepackage{amsfonts,amsthm}
\usepackage{array}
\usepackage{amscd}
\usepackage{amssymb}
\usepackage{enumerate}

\newtheorem{theorem}{Theorem}
\newtheorem{lemma}{Lemma}
\newtheorem{corollary}{Corollary}
\newtheorem{remark}{Remark}

\newcounter{ex}
\newenvironment{example}{\par\noindent\textbf{Example \stepcounter{ex}\arabic{ex}.}}{\bigskip}

\input{tcilatex}
\begin{document}

\title{On a global implicit function theorem for locally Lipschitz maps via
nonsmooth critical point theory}
\author{M. Galewski and M. R\u{a}dulescu}
\maketitle

\begin{abstract}
\noindent We prove a non-smooth generalization of the global implicit
function theorem. More precisely we use the non-smooth local implicit
function theorem and the non-smooth critical point theory in order to prove
a non-smooth global implicit function theorem for locally Lipschitz
functions. A comparison between several global inversion theorems is
discussed.
\end{abstract}

\section{Introduction}

In this note we provide conditions for the existence of a global implicit
function for the equation $F\left( x,y\right) =0$, where $F:\mathbb{R}%
^{n}\times \mathbb{R}^{m}\rightarrow \mathbb{R}^{m}$ is a locally Lipschitz
functions. We use non-smooth critical point theory together with classical
tools for the existence of a local (non-smooth) implicit function. The
following theorem is a finite dimensional counterpart of the main result
given in Idczak \cite{idczak}.

\begin{theorem}
\label{IdczakIFTcit}Assume that $F:$ $\mathbb{R}^{n}\times \mathbb{R}%
^{m}\rightarrow \mathbb{R}^{n}$ is a $C^{1}$ mapping such that:\medskip 
\newline
(a1) for any $y\in \mathbb{R}^{m}$ the functional $\varphi _{y}:\mathbb{R}%
^{n}\rightarrow \mathbb{R}$ given by the formula 
\begin{equation*}
\varphi _{y}\left( x\right) =\frac{1}{2}\left\Vert F\left( x,y\right)
\right\Vert ^{2}
\end{equation*}%
is coercive, i.e. $\lim_{\left\Vert x\right\Vert \rightarrow \infty }\varphi
_{y}\left( x\right) =+\infty $,\medskip \newline
(a2) the Jacobian matrix $F_{x}(x,y)$ is bijective for any $(x,y)\in \mathbb{%
R}^{n}\times \mathbb{R}^{m}$.\medskip \newline
Then there exists a unique function $f:\mathbb{R}^{m}\rightarrow \mathbb{R}%
^{n}$ such that equations $F(x,y)=0$ and $x=f(y)$ are equivalent in the set $%
\mathbb{R}^{n}\times \mathbb{R}^{m}$, in other words  $F(f\left( y\right)
,y)=0$ for any $y\in \mathbb{R}^{m}$. Moreover, $f\in C^{1}(\mathbb{R}^{m},%
\mathbb{R}^{n})$.
\end{theorem}

The proof of Theorem \ref{IdczakIFTcit} is based on the application of a
local implicit function theorem in the $C^{1}$ setting and next on the
application of classical mountain pass theorem. In this paper we intend to
provide a direct proof of its locally Lipschitz version. Thus we coin
together local implicit function results and the non-smooth mountain pass
technique. Obtaining a non-smooth version of the infinite dimensional
counterpart of result in \cite{idczak} is an open problem. The reason is
that the non-smooth results are based on the Clarke's local inversion
theorem for locally Lipschitz functions, which works only in a finite
dimensional framework.\bigskip 

In order to prove our results we start with some preliminaries.\bigskip 

A function $f:%
\mathbb{R}
^{n}\rightarrow 
\mathbb{R}
^{k}$ is called locally Lipschitz continuous, if to every $u\in 
\mathbb{R}
^{n}$ there corresponds a neighborhood $V_{u}$ of $u$ and a constant $%
L_{u}\geq 0$ such that 
\begin{equation*}
\left\Vert f(z)-f(w)\right\Vert \leq L_{u}\Vert z-w\Vert \quad \text{%
for\thinspace all \ }z,w\in V_{u}\,.
\end{equation*}%
If $u,z\in 
\mathbb{R}
^{n}$, we write $f^{0}(u;z)$ for the generalized directional derivative of $%
f:%
\mathbb{R}
^{n}\rightarrow 
\mathbb{R}
$ at the point $u$ along the direction $z$; i.e., 
\begin{equation*}
f^{0}(u;z):=\limsup_{w\rightarrow u,\,t\rightarrow 0^{+}}\frac{f(w+tz)-f(w)}{%
t}\,.
\end{equation*}%
The generalized gradient of the function $f$ at $u$, denoted by $\partial
f(u)$, is the set 
\begin{equation*}
\partial f(u):=\{\xi \in L\left( 
\mathbb{R}
^{n},%
\mathbb{R}
\right) :\langle \xi ,z\rangle \leq f^{0}(u;z),\;\text{for all }\,z\in 
\mathbb{R}
^{n}\}.
\end{equation*}%
For the definition of a generalized Jacobian of a vector valued function $f:%
\mathbb{R}
^{n}\rightarrow 
\mathbb{R}
^{k}$\ we refer to \cite{clarkebook} p. 69. We denote the generalized
Jacobian at $x$\ again by $\partial f(x)$. For a fixed $x$\ the set $%
\partial f(x)$\ being of maximal rank means that all matrices in $\partial
f(x)$\ are nonsingular.\texttt{\ }This assumption is equivalent, when $f$ is
smooth, with the assumption that $\det [f^{\prime }(x)]\neq 0$ for every $%
x\in D$ where $D\subset 
\mathbb{R}
^{n}$\ is some open set. Compare with \cite{rad} where this condition
provides local diffeomorphism for a differentiable mapping. Note that it is
not enough to assume that $\det [f^{\prime }(x)]\neq 0$ whenever it exists,
which happens a.e. for a locally Lipschitz function.\medskip 

The basic properties of generalized directional derivative and generalized
gradient were studied in \cite{clarkebook} and later in \cite{ochalbook}.
Now we provide the local implicit function in the sense of Clarke as
suggested in \cite{sun}.

\begin{theorem}
\label{finite_dimensional_implicit copy(1)} \cite{sun}Assume that $F:\mathbb{%
R}^{n}\times \mathbb{R}^{m}\rightarrow \mathbb{R}^{n}$ is a locally
Lipschitz mapping in a neigbourhood of a point $\left( x_{0},y_{0}\right) $
such that $F\left( x_{0},y_{0}\right) =0$. Assume further that $\partial
_{x}F(x_{0},y_{0})$ is of maximal rank. Then there exist a neighborhood $%
V\subset \mathbb{R}^{m}$ of $y_{0}$ and a Lipschitz function $G:V\rightarrow 
\mathbb{R}^{n}$ such that for every $y$ in $V$ it holds $F(G(y),y)=0$ and $%
G(y_{0})=x_{0}$.
\end{theorem}

We will later on make this result a global one with additional coercivity
assumption.

\section{Global implicit function result}

A point $u$ is called a (generalized) critical point of the locally
Lipschitz continuous functional $J:%
\mathbb{R}
^{n}\rightarrow 
\mathbb{R}
$ if $0\in \partial J(u)$; in this case we identify $L\left( 
\mathbb{R}
^{n},%
\mathbb{R}
\right) $ with $%
\mathbb{R}
^{n}$ so that $\partial J(u)\subset 
\mathbb{R}
^{n}$. $J$ is said to fulfill the non-smooth Palais-Smale condition, see 
\cite{MoVa}, if every sequence $\{u_{n}\}$ in $%
\mathbb{R}
^{n}$ such that $\{J(u_{n})\}$ is bounded and 
\begin{equation*}
J^{0}(u_{n};u-u_{n})\geq -\varepsilon _{n}\Vert u-u_{n}\Vert 
\end{equation*}%
\ for all $u\in 
\mathbb{R}
^{n}$, where $\varepsilon _{n}\rightarrow 0^{+}$, admits a convergent
subsequence. Our main tool will be the following result based on the
zero-altitude version of Mountain Pass Theorem from \cite{MoVa}, where we
replace non-smooth (PS) condition with coercivity which we require and which
guarantees that (PS) condition holds.

\begin{theorem}
\label{MPT} Let $J:%
\mathbb{R}
^{n}\rightarrow \mathbb{R}$ be a coercive locally Lipschitz continuous
functional. If there exist $u_{1},u_{2}\in 
\mathbb{R}
^{n}$, $u_{1}\neq u_{2}$ and $r\in (0,\left\Vert u_{2}-u_{1}\right\Vert )$
such that 
\begin{equation*}
\inf \{J(u):\left\Vert u-u_{1}\right\Vert =r\}\geq \max \{J(u_{1}),J(u_{2})\}
\end{equation*}%
and we denote by $\Gamma $ the family of continuous paths $\gamma
:[0,1]\rightarrow 
\mathbb{R}
^{n}$ joining $u_{1}$ and $u_{2},$ then%
\begin{equation*}
c:=\underset{\gamma \in \Gamma }{\inf }\underset{s\in \lbrack 0,1]}{\max }%
J(\gamma (s))\geq \max \{J(u_{1}),J(u_{2})\}
\end{equation*}%
is a critical value for $J$ on $%
\mathbb{R}
^{n}$\ and $K_{c}\backslash \{u_{1},u_{2}\}\neq \emptyset $, where $K_{c}$
is the set of critical points at the level $c$, i.e.%
\begin{equation*}
K_{c}=\left\{ u\in 
\mathbb{R}
^{n}:J\left( u\right) =c\text{ and }0\in \partial J(u)\right\} .
\end{equation*}
\end{theorem}

We are now able to formulate and prove our main result which is a non-smooth
generalization of a global implicit function theorem from \cite{chua}.

\begin{theorem}
\label{finite_dimensional_implicit} Assume that $F:\mathbb{R}^{n}\times 
\mathbb{R}^{m}\rightarrow \mathbb{R}^{n}$ is a locally Lipschitz mapping
such that\newline
b1) For every $y\in \mathbb{R}^{m}$ the functional $\varphi _{y}:\mathbb{R}%
^{n}\rightarrow \mathbb{R}$ given by the formula 
\begin{equation*}
\varphi _{y}\left( x\right) =\frac{1}{2}\left\Vert F\left( x,y\right)
\right\Vert ^{2}
\end{equation*}%
is coercive, i.e. $\lim_{\left\Vert x\right\Vert \rightarrow \infty }\varphi
_{y}\left( x\right) =+\infty $,\newline
b2) the set $\partial _{x}F(x,y)$ is of maximal rank for all $(x,y)\in 
\mathbb{R}^{n}\times \mathbb{R}^{m}$.\newline
Then there exists a unique locally Lipschitz function $f:\mathbb{R}%
^{m}\rightarrow \mathbb{R}^{n}$ such that equations $F(x,y)=0$ and $x=f(y)$
are equivalent in the set $\mathbb{R}^{n}\times \mathbb{R}^{m}$.
\end{theorem}

\textbf{Proof. }We follow the ideas used in the proof of Main Theorem in 
\cite{idczak} with necessary modifications due to the fact that we now use
non-smooth critical point theory. In view of Theorem \ref%
{finite_dimensional_implicit copy(1)} assumption (\textit{b2}) implies that $%
F\left( \cdot ,y\right) $ defines a locally invertible mapping with any
fixed $y\in 
\mathbb{R}
^{m}$. Thus it is sufficient to show that this map is onto and one to one.
\bigskip

Let us fix a point $y\in 
\mathbb{R}
^{m}$. Functional $\varphi _{y}$ is continuous and coercive, so it has an
argument of a minimum $\overline{x}$ which necessarily is the (generalized)
critical point, i.e. $0\in \partial \varphi _{y}\left( \overline{x}\right) $%
. By Lemma 2 from \cite{CLARKE2} (or else Theorem 2.6.6. p. 85 from \cite%
{clarkebook})\texttt{\ }noting that a $C^{1}$\ function is strictly
differentiable, we obtain 
\begin{equation*}
0\in \partial _{x}\left( \frac{1}{2}\left\Vert F\left( \overline{x},y\right)
\right\Vert ^{2}\right) \subset \partial _{x}F(\overline{x},y)\circ F\left( 
\overline{x},y\right) \text{.}
\end{equation*}%
Since by assumption b2) set $\partial _{x}F(\overline{x},y)$ is of maximal
rank, we see\texttt{\ }that it must hold that\texttt{\ }$F\left( \overline{x}%
,y\right) =0$ and so to each $y$ there corresponds at least one $x$ such
that $F\left( x,y\right) =0$ is solvable. \bigskip

Now we argue by contradiction in order to show that $F\left( \cdot ,y\right) 
$ is one to one. Suppose that there are $x_{1},x_{2}\in \mathbb{R}^{n}$, $%
x_{1}\neq x_{2}$, such that $F(x_{1},y)=F(x_{2},y)=0$. We will apply Theorem %
\ref{MPT}. Thus we put $e=x_{1}-x_{2}$ and define mapping $g_{y}:%
\mathbb{R}
^{n}\rightarrow 
\mathbb{R}
^{k}$ by the following formula%
\begin{equation*}
g_{y}\left( x\right) =F\left( x+x_{2},y\right)
\end{equation*}%
Note that $g_{y}\left( 0\right) =g_{y}\left( e\right) =0$. We define a
locally Lipschitz functional $\psi _{y}:%
\mathbb{R}
^{n}\rightarrow 
\mathbb{R}
$ by 
\begin{equation*}
\psi _{y}\left( x\right) =\frac{1}{2}\left\Vert g_{y}\left( x\right)
\right\Vert ^{2}.
\end{equation*}%
By assumption (\textit{b1}) and by its definition functional $\psi _{y}$ is
coercive and also $\psi _{y}\left( e\right) =\psi _{y}\left( 0\right) =0.$%
\texttt{\ Fix }$\rho >0$\texttt{\ }such that $\rho <\left\Vert e\right\Vert $%
\ and consider functional $\psi _{y}$\ on the boundary of $\overline{B\left(
0,\rho \right) }$. Since $\psi _{y}$\ is continuous, it attains its minimum
at a point $x=\overline{x}$. We claim that there exists at least one such
radius $\rho $\ that $\inf_{\left\Vert x\right\Vert =\rho }\psi _{y}(x)>0$.
Suppose this is not the case and that for any $\rho <\left\Vert e\right\Vert 
$\ we have $\inf_{\left\Vert x\right\Vert =\rho }\psi _{y}(x)=0$. This means
there exists $\overline{x}$\ from the boundary of $\overline{B\left( 0,\rho
\right) }$\ such that $\psi _{y}(\overline{x})=0$. Then $F\left( \overline{x}%
+x_{2},y\right) =0$\ and since $y$\ is held fixed we obtain contradiction
with the local implicit function theorem, namely Theorem \ref%
{finite_dimensional_implicit copy(1)}. Thus there exists some $\rho >0$\
such that\texttt{\ }%
\begin{equation*}
\inf_{\left\Vert x\right\Vert =\rho }\psi _{y}(x)>0=\psi _{y}\left( e\right)
=\psi _{y}\left( 0\right)
\end{equation*}%
Thus by Theorem \ref{MPT} applied to $J=\psi _{y}$\ we note that $\psi _{y}$%
\ has a generalized critical point $v$\ which is different from $0$\ and $e$%
\ since the corresponding critical value $\psi _{y}\left( v\right) \geq
\inf_{\left\Vert x\right\Vert =\rho }\psi _{y}(x)>0$. On the other hand,
since $v$\ is a critical point, we get by the chain rule from \cite{CLARKE2}
mentioned at the beginning of the proof that 
\begin{equation*}
0\in \partial \psi _{y}(v)\subset \partial _{x}F(v+x_{2},y)\circ F\left(
v+x_{2},y\right) \text{.}
\end{equation*}%
Since $\partial _{x}F(v+x_{2},y)$\ is of full rank, we see that $F\left(
v+x_{2},y\right) =0$. This means that the equality $\psi _{y}(v)=0$\ holds
which contradicts $\psi _{y}(v)>0$. The obtained contradiction ends the
proof. \hfill $\blacksquare $

\section{A remark on a non-smooth generalization of the Hadamard-Palais
theorem}

The following theorem is known as the Hadamard-Palais theorem

\begin{theorem}
\label{MainTheo copy(2)}\cite{palais}Let $X,$ $B$ be finite dimensional
Banach spaces. Assume that $f:X\rightarrow B$ is a $C^{1}$-mapping such that%
\newline
(c1) $\left\Vert f\left( x\right) \right\Vert \rightarrow \infty $ as $%
\left\Vert x\right\Vert \rightarrow \infty $;\newline
(c2) $f^{\prime }(x)$ is invertible for any $x\in X$\newline
then $f$ is a diffeomorphism.
\end{theorem}

With our approach we are able to provide a direct proof of the locally
Lipschitz version of the above theorem

\begin{theorem}
\label{MainTheo} Assume $f:%
\mathbb{R}
^{n}\rightarrow 
\mathbb{R}
^{n}$ is a locally Lipschitz mapping such that\newline
(d1) for any $y\in 
\mathbb{R}
^{n}$ the functional $\varphi _{y}:%
\mathbb{R}
^{n}\rightarrow 
\mathbb{R}
$ defined by 
\begin{equation*}
\varphi _{y}\left( x\right) =\frac{1}{2}\left\Vert f\left( x\right)
-y\right\Vert ^{2}
\end{equation*}%
is coercive;\newline
(d2) for any $x\in 
\mathbb{R}
^{n}$ we have that $\partial f(x)$ is of maximal rank.\newline
Then $f$ is a global homeomorphism on $%
\mathbb{R}
^{n}$ and $f^{-1}$ is locally Lipschitz.
\end{theorem}

However, one remark is in order. Let us recall the following

\begin{lemma}
\cite{CLARKE2},\cite{clarkebook}\label{lemma2}Let $D$ be an open subset of $%
\mathbb{R}
^{n}$. If $f:D\rightarrow 
\mathbb{R}
^{n}$ satisfies a Lipschitz condition in some neighbourhood of $x_{0}$ and $%
\partial f(x_{0})$ is of maximal rank, then there exist neighborhoods $%
U\subset D$ of $x_{0}$ and $V$ of $f(x_{0})$ and a Lipschitz function $%
g:V\rightarrow 
\mathbb{R}
^{n}$ such that\newline
i ) for every $u\in U$, $g(f(u))=u$, and\newline
ii) for every $v\in V$, $f(g(v))=v.$
\end{lemma}

Now, note that condition (d1) is equivalent to (c1) and condition (d2)
implies $f$ is a local homeomorphism. By the Banach-Mazur Theorem it then
follows that $f$ is a global homeomorphism, so it is well known that $f^{-1}$
is continuous. From Lemma \ref{lemma2} it then follows that $f^{-1}$ is also
locally Lipschitz.\bigskip 

\subsection{Applications to algebraic equations}

We finish this section with some applications of Theorem \ref{MainTheo} to
algebraic equations which received some interest in the literature lately,
see for example \cite{add3}, \cite{add4}, \cite{add2}, \cite{add1} where
various variational approaches are being applied to the existence and
multiplicity.\bigskip 

In the following we consider the problem%
\begin{equation}
Ax=F\left( x\right) +\xi ,  \label{equ_App_1}
\end{equation}%
where $\xi \in 
\mathbb{R}
^{N}$ is fixed, $A$ is an $N\times N$ matrix which is not positive definite,
negative definite or symmetric; $F:%
\mathbb{R}
^{N}\rightarrow 
\mathbb{R}
^{N}$ is a locally Lipschitz function. We consider $%
\mathbb{R}
^{N}$ with Euclidean norm in both theoretical results and the example which
follows.

Note that when A is such as above one cannot apply even the simplest
variational approach, i.e. the direct method relying on minimizing the Euler
action functional%
\begin{equation*}
J\left( x\right) =\left\langle Ax,x\right\rangle -\mathcal{F}\left( x\right)
,
\end{equation*}%
and where $\mathcal{F}:%
\mathbb{R}
^{N}\rightarrow 
\mathbb{R}
$ is the potential of $F$. The difficulties are due to the fact that term $%
\left\langle Ax,x\right\rangle $ need not be coercive nor anti-coercive.
Moreover, uniqueness which we achieve, in the classical approach requires
strict convexity of the action functional which is again an assumption
rather demanding.\bigskip 

In order apply Theorem \ref{MainTheo} to the solvability of (\ref{equ_App_1}%
) we need some assumptions. Let us recall that if $A^{\ast }$ denotes the
transpose of matrix $A$, then $A^{\ast }A$ is symmetric and positive
semidefinite. However, $A^{\ast }A$ being positive semidefinite is not
sufficient for our purposes. We assume what follows\bigskip

\begin{description}
\item[A1] Matrix $A^{\ast }A$ is positive definite with eigenvalues ordered
as 
\begin{equation*}
0<\lambda _{1}\leq \cdot \cdot \cdot \leq \lambda _{N}.
\end{equation*}%
\bigskip 
\end{description}

Now we can state the following existence theorems.

\begin{theorem}
\label{firstAlgTheo}Assume that \textbf{A1} holds, $F:%
\mathbb{R}
^{N}\rightarrow 
\mathbb{R}
^{N}$ is a locally Lipschitz function and that the following conditions hold%
\newline
(i) There exists a constant $0<a<\sqrt{\lambda _{1}}$ such that%
\begin{equation*}
\left\Vert F\left( x\right) \right\Vert \leq a\left\Vert x\right\Vert 
\end{equation*}
for all sufficiently large $x\in 
\mathbb{R}
^{N},$\newline
(ii) $\det (A-B)\neq 0$ for every $B\in \partial F\left( x\right) $ and
every $x\in 
\mathbb{R}
^{N}.$\newline
Then problem (\ref{equ_App_1}) has exactly one solution  for any $\xi \in 
\mathbb{R}
^{N}$.
\end{theorem}

\textbf{Proof.  }We need to show that assumptions of Theorem \ref{MainTheo}
are satisfied. We put $\varphi \left( x\right) =Ax-F\left( x\right) $. In
order to demonstrate (d1) we see that for sufficiently large $x\in 
\mathbb{R}
^{N}$ 
\begin{equation*}
\begin{array}{l}
\left\Vert \varphi \left( x\right) \right\Vert =\left\Vert Ax-F\left(
x\right) \right\Vert \geq \left\Vert Ax\right\Vert -\left\Vert F\left(
x\right) \right\Vert \geq \bigskip  \\ 
\sqrt{\left\langle A^{\ast }Ax,x\right\rangle }-a\left\Vert x\right\Vert
\geq \left( \sqrt{\lambda _{1}}-a\right) \left\Vert x\right\Vert .%
\end{array}%
\end{equation*}%
\hfill Hence the function $\varphi $ is coercive. From (ii) it follows that
the $\partial \varphi \left( x\right) $ is of maximal rank for all $x\in 
\mathbb{R}
^{N}$. Therefore condition (d2) is satisfied. From Theorem \ref{MainTheo} it
follows that $\varphi $ is a global homeomorphism and equation (\ref%
{equ_App_1}) has exactly one solution for any $\xi \in 
\mathbb{R}
^{N}$. $\blacksquare $\textbf{\ } 

\begin{theorem}
\label{SecondAlegTheo}Assume that \textbf{A1} holds, $F:%
\mathbb{R}
^{N}\rightarrow 
\mathbb{R}
^{N}$ is a locally Lipschitz function and that the following conditions hold%
\newline
(i) There exists a constant $b>\sqrt{\lambda _{N}}$ such that%
\begin{equation*}
\left\Vert F\left( x\right) \right\Vert \geq b\left\Vert x\right\Vert 
\end{equation*}%
for all sufficiently large $x\in 
\mathbb{R}
^{N},$\newline
(ii) $\det (A-B)\neq 0$ for every $B\in \partial F\left( x\right) $ and
every $x\in 
\mathbb{R}
^{N}.$\newline
Then problem (\ref{equ_App_1}) has exactly one solution for any $\xi \in 
\mathbb{R}
^{N}$.
\end{theorem}

\textbf{Proof. }We put $\varphi _{1}\left( x\right) =F\left( x\right) -Ax$
and we observe that for sufficiently large $x\in 
\mathbb{R}
^{N}$ 
\begin{equation*}
\left\Vert \varphi _{1}\left( x\right) \right\Vert =\left\Vert F\left(
x\right) -Ax\right\Vert \geq \left\Vert F\left( x\right) \right\Vert
-\left\Vert Ax\right\Vert \geq \left( b-\sqrt{\lambda _{N}}\right)
\left\Vert x\right\Vert .
\end{equation*}%
\hfill Hence the function $\varphi _{1}$ is coercive and the assertion
follows as in the proof of the above result.  $\blacksquare $\textbf{\ } 

\begin{remark}
We note that in order to get coercivity of function $\varphi $ in Theorem %
\ref{firstAlgTheo} we can replace condition (i) with the following
assumption: exist constants $\alpha >0$, $0<\gamma <1$ such that%
\begin{equation*}
\left\Vert F\left( x\right) \right\Vert \leq \alpha \left\Vert x\right\Vert
^{\gamma }
\end{equation*}
for all sufficiently large $x\in 
\mathbb{R}
^{N}$. Indeed, we have for sufficiently large $x\in 
\mathbb{R}
^{N}$%
\begin{equation*}
\left\Vert \varphi \left( x\right) \right\Vert \geq \left\Vert x\right\Vert
^{\gamma }\left( \sqrt{\lambda _{1}}\left\Vert x\right\Vert ^{1-\gamma
}-\alpha \right) ,
\end{equation*}%
which means that $\varphi $ is coercive. \bigskip \newline
Concerning Theorem \ref{SecondAlegTheo} we can replace condition (i) with
the following assumption: there exist constants $\beta >0$, $\theta >1$ such
that%
\begin{equation*}
\left\Vert F\left( x\right) \right\Vert \leq \beta \left\Vert x\right\Vert
^{\theta }
\end{equation*}
for all sufficiently large $x\in 
\mathbb{R}
^{N}$. In this case we obtain%
\begin{equation*}
\left\Vert \varphi _{1}\left( x\right) \right\Vert \geq \left( \beta
\left\Vert x\right\Vert ^{\theta -1}-\sqrt{\lambda _{N}}\right) \left\Vert
x\right\Vert 
\end{equation*}%
for all sufficiently large $x\in 
\mathbb{R}
^{N}$, which means that $\varphi _{1}$ is coercive. 
\end{remark}

\begin{example}
Consider an indefinite matrix $A=\left[ 
\begin{array}{ll}
-2 & 1 \\ 
4 & -3%
\end{array}%
\right] $ and function $F:%
\mathbb{R}
^{2}\rightarrow 
\mathbb{R}
^{2}$ given by 
\begin{equation*}
F\left( x,y\right) =\left( x^{3}+\left\vert y\right\vert ,4x+\left\vert
y\right\vert +y^{3}\right) 
\end{equation*}%
Consider on $%
\mathbb{R}
^{2}$ the Euclidean norm that is $\left\Vert \left( x,y\right) \right\Vert =%
\sqrt{x^{2}+y^{2}}$ and we recall that 
\begin{equation*}
\left\Vert \left( x,y\right) \right\Vert \leq 2^{\frac{1}{3}}\sqrt[6]{%
x^{6}+y^{6}}
\end{equation*}%
We use also symbol $\mathbf{u=}\left( x,y\right) $. Note that for any $%
\mathbf{u}\in 
\mathbb{R}
^{2}$%
\begin{equation*}
\begin{array}{l}
\left\Vert \left( x^{3}+\left\vert y\right\vert ,4x+\left\vert y\right\vert
+y^{3}\right) \right\Vert \geq \left\Vert \left( x^{3},4x+y^{3}\right)
\right\Vert -\left\Vert \left( \left\vert y\right\vert ,\left\vert
y\right\vert \right) \right\Vert \geq \bigskip  \\ 
\sqrt{x^{6}+y^{6}}-2\left\vert y\right\vert \geq \sqrt{x^{6}+y^{6}}%
-2\left\Vert \mathbf{u}\right\Vert \geq \frac{1}{2}\left\Vert \mathbf{u}%
\right\Vert ^{3}-2\left\Vert \mathbf{u}\right\Vert 
\end{array}%
\end{equation*}%
Moreover, we can directly calculate that $\left\Vert A\mathbf{u}\right\Vert $
$\leq \left\Vert A\right\Vert \left\Vert \mathbf{u}\right\Vert $. This means
that 
\begin{equation*}
\left\Vert \left( x^{3}+\left\vert x\right\vert ,4x+\left\vert x\right\vert
+y^{3}\right) -A\mathbf{u}\right\Vert \geq \frac{1}{2}\left\Vert \mathbf{u}%
\right\Vert ^{3}-\left\Vert A\right\Vert \left\Vert \mathbf{u}\right\Vert
=\left( \frac{1}{2}\left\Vert \mathbf{u}\right\Vert ^{2}-\left\Vert
A\right\Vert \right) \left\Vert \mathbf{u}\right\Vert .
\end{equation*}%
Hence for sufficiently large $\left\Vert \mathbf{u}\right\Vert $, we see
that 
\begin{equation*}
\left\Vert \left( x^{3}+\left\vert x\right\vert ,4x+\left\vert x\right\vert
+y^{3}\right) -A\mathbf{u}\right\Vert \geq \frac{1}{3}\left\Vert \mathbf{u}%
\right\Vert ^{3}
\end{equation*}%
and so assumption \textbf{A4} holds. As for assumption \textbf{A3} we see
that $B\in \partial F\left( \mathbf{u}\right) $ has the following form%
\begin{equation*}
B=\left[ 
\begin{array}{cc}
3x^{2} & 0 \\ 
4 & 3y^{2}%
\end{array}%
\right] +\left[ 
\begin{array}{ll}
0 & a \\ 
a & 0%
\end{array}%
\right] ,
\end{equation*}%
where $a\in \left[ -1,1\right] $ for any $\mathbf{u}\in 
\mathbb{R}
^{2}$. Note also that 
\begin{equation*}
B-A=\left[ 
\begin{array}{cc}
3x^{2}+2 & a-1 \\ 
a & 3y^{2}+3%
\end{array}%
\right] .
\end{equation*}%
We calculate that $\det \left( B-A\right) =\left( 3x^{2}+2\right) \left(
3y^{2}+3\right) -a\left( a-1\right) $. Since $a\in \left[ -1,1\right] $, we
see $\det \left( B-A\right) >0$. \bigskip 
\end{example}

Sometimes it is easier to prove coercivity of $\left\Vert Ax-F\left(
x\right) \right\Vert $ directly than to use growth conditions on the
nonlinear term. Moreover, when we prove the coercivity directly, there is no
need to assume that $A^{\ast }A$ is positive definite. Thus from the proof
of Theorems \ref{firstAlgTheo} and \ref{SecondAlegTheo} it follows that

\begin{corollary}
\label{corollAlgebEqu}Assume that\newline
(i) $\left\Vert Ax-F\left( x\right) \right\Vert \rightarrow \infty $ as $%
\left\Vert x\right\Vert \rightarrow \infty ,$ \newline
(ii) $\det (A-B)\neq 0$ for every $B\in \partial F\left( x\right) $ and
every $x\in 
\mathbb{R}
^{N}.$\newline
Then (\ref{equ_App_1}) has exactly one solution for any fixed $\xi \in 
\mathbb{R}
^{N}$.
\end{corollary}

Now we proceed with providing example of a nonlinear term $F$ and a matrix $%
A $ for which our above result works.

\begin{example}
Consider an indefinite matrix $A=\left[ 
\begin{array}{cc}
-\frac{3}{2} & 4 \\ 
5 & -\frac{40}{3}%
\end{array}%
\right] $ such that $\det A=0$ and a function $F:%
\mathbb{R}
^{2}\rightarrow 
\mathbb{R}
^{2}$ given by 
\begin{equation*}
F\left( x,y\right) =\left( x^{5}+\left\vert y\right\vert ,4x+\left\vert
y\right\vert +y^{5}\right) 
\end{equation*}%
Note that for any $\mathbf{u}\in 
\mathbb{R}
^{2}$%
\begin{equation*}
\begin{array}{l}
\left\Vert \left( x^{5}+\left\vert y\right\vert ,4x+\left\vert y\right\vert
+y^{5}\right) \right\Vert \geq \left\Vert \left( x^{5},y^{5}+4x\right)
\right\Vert -\left\Vert \left( \left\vert y\right\vert ,\left\vert
y\right\vert \right) \right\Vert \geq \bigskip  \\ 
\left\Vert \left( x^{5},y^{5}\right) \right\Vert -\left\Vert \left(
0,4x\right) \right\Vert -\left\Vert \left( \left\vert y\right\vert
,\left\vert y\right\vert \right) \right\Vert =\sqrt{x^{10}+y^{10}}-\left( 4+%
\sqrt{2}\right) \left\vert y\right\vert 
\end{array}%
\end{equation*}%
Moreover%
\begin{equation*}
\left\Vert \left[ 
\begin{array}{cc}
-\frac{3}{2} & 4 \\ 
5 & -\frac{40}{3}%
\end{array}%
\right] \left[ 
\begin{array}{l}
x \\ 
y%
\end{array}%
\right] \right\Vert \leq \frac{13}{2}\left\vert x\right\vert +\frac{52}{3}%
\left\vert y\right\vert .
\end{equation*}%
Thus we have 
\begin{equation*}
\left\Vert A\mathbf{u}-F\left( \mathbf{u}\right) \right\Vert \geq \sqrt{%
x^{10}+y^{10}}-\left( 4+\sqrt{2}-\frac{52}{3}\right) \left\vert y\right\vert
-\frac{13}{2}\left\vert x\right\vert \rightarrow \infty 
\end{equation*}%
as $\left\Vert \left( x,y\right) \right\Vert \rightarrow \infty $.\bigskip
Suppose this is not the case. Then there is some $M>0$ such that for any $R>0
$ and all $\left\Vert \left( x,y\right) \right\Vert \geq R$ it holds $\sqrt{%
x^{10}+y^{10}}-\left( 4+\sqrt{2}-\frac{52}{3}\right) \left\vert y\right\vert
-\frac{13}{2}\left\vert x\right\vert <M$. So we easily arrive at a
contradiction.\bigskip \newline
\qquad Observe that $B\in \partial F\left( \mathbf{u}\right) $ has the
following form%
\begin{equation*}
B=\left[ 
\begin{array}{ll}
5x^{4} & 0 \\ 
4 & 5y^{4}%
\end{array}%
\right] +\left[ 
\begin{array}{ll}
0 & a \\ 
a & 0%
\end{array}%
\right] ,
\end{equation*}%
where $a\in \left[ -1,1\right] $ for any $\mathbf{u}\in 
\mathbb{R}
^{2}$. Note also that 
\begin{equation*}
B-A=\left[ 
\begin{array}{cc}
5x^{4}+\frac{3}{2} & a-4 \\ 
a-1 & 5y^{4}+\frac{40}{3}%
\end{array}%
\right] .
\end{equation*}%
Since a function $f\left( x\right) =\left( x-1\right) \left( x-4\right) $
has a maximal value $10$ on $\left[ -1,1\right] $ and since $\left( 5x^{4}+%
\frac{3}{2}\right) \left( y^{4}+\frac{40}{3}\right) \geq 0$, we see that $%
\det \left( B-A\right) \neq 0$. Note that in this case $A^{\ast }A$ does not
satisfy \textbf{A1}.
\end{example}

\section{A comparison between several global inversion theorems}

It is well-known that the implicit function and inverse theorems are
equivalent in the sense that the validity of one implies the validity of the
other. In the following we shall show the novelty of our main result. We
shall make comparisons between several global inversion theorems known from
the literature.

The following theorem is known as the Hadamard-Levy theorem (see \cite%
{Hadamard}, \cite{Levy}, \cite{plastock}, \cite{rad2})

\begin{theorem}[\textbf{Hadamard-Levy theorem}]
Let $E,F$ be two Banach spaces and $f\ :E\rightarrow F$ be \textsl{.}a local
diffeomorphism of class $C^{1}$ which satisfies the following integral
condition 
\begin{equation*}
\int_{0}^{\infty }\min_{\left\Vert x\right\Vert =r}\Vert f\ ^{\prime
}(x)^{-1}\Vert ^{-1}dr=\infty .
\end{equation*}%
Then $f$\ is a global diffeomorphism.
\end{theorem}

The finite dimensional version of the above theorem was extended to locally
Lipschitz functions by Pourciau, see \cite{pourcianu1}, \cite{pourcianu2}.
If $A$ is a square matrix we denote $\left[ A\right] =\underset{\left\Vert
u\right\Vert =1}{\inf }\left\Vert Au\right\Vert .$

\begin{theorem}[\textbf{Pourciau's theorem}]
Let $f:%
\mathbb{R}
^{n}\rightarrow 
\mathbb{R}
^{n}$ be a locally Lipschitz function and suppose that the generalized
Jacobian $\partial f\left( x\right) $ is of full rank for every $x\in 
\mathbb{R}
^{n}$. Let $m\left( t\right) =\underset{\left\Vert z\right\Vert \leq \ t}{%
\inf }\left[ \partial f\left( z\right) \right] =\underset{\left\Vert
z\right\Vert \leq \ t}{\inf }\underset{\ \ \ A\in \partial f\left( z\right) }%
{\inf }\left[ A\right] $ and suppose that%
\begin{equation*}
\overset{\infty }{\underset{0}{\dint }}m\left( t\right) dt=+\infty
\end{equation*}%
Then $f$ is a bijective function and the inverse of $f$, that is $f^{-1}$,
is a locally Lipschitz function.
\end{theorem}

One interesting global invertibility result for nonsmooth functions was
stated in \cite{ioffee}. In order to state the result we shall define the
modulus of surjection of a function $f$ at a point $x.$ Let $E,F$ be two
Banach spaces, $f:E\rightarrow F$ and $x\in E$. We denote by $B[a,r]$ the
closed ball of radius $r$ around $a\in E.$

\begin{equation*}
\text{Sur}(f,x)(t)=\sup \{r\geq 0:B[f(x),r]\subset f[B(x,t)]\}>o
\end{equation*}

Thus, for any $t$ $>0$, the value of the modulus of surjection of $f$ at $x$
is the maximal radius of a ball around $f(x)$ contained in the $f$-image of
the ball of radius $t$ around $x$. We further introduce the constant of
surjection of $f$ at $x$ by

\begin{equation*}
\text{sur}\left( f,x\right) =\underset{t\rightarrow 0}{\lim \inf }\frac{%
\text{Sur}(f,x)(t)}{t}
\end{equation*}

Obviously, sur$\left( f,x\right) >0$ is a sufficient condition for $f$ to be
surjective at $x$, that is, for Sur$(f,x)(t)$ to be positive for small $t.$
The following two theorems are taken from \cite{ioffee}.

\begin{theorem}
Let $E,F$ be two Banach spaces and $f:E\rightarrow F.$ Suppose the graph of $%
f$ is closed and there is a positive lower semicontinuous (l.s.c.) function $%
m:[0,\infty )\rightarrow \lbrack 0,\infty )$ such that%
\begin{equation}
\text{sur}\left( f,x\right) \geq m\left( \left\Vert x\right\Vert \right) ,\
x\in E\text{ }  \label{surf(x,t)3.1}
\end{equation}%
Then
\end{theorem}

\begin{equation*}
\text{Sur}(f,x)(t)\geq \underset{0}{\overset{r}{\dint }}m\left( s\right) ds\
\ \ \ \ \ \ \text{for every }r>0
\end{equation*}

\begin{theorem}
\label{TheoFromIoffe}Let $E,F$ be two Banach spaces and $f:E\rightarrow F$
be a continuous mapping that is locally one-to-one (i.e., every $x\in E$ has
a neighborhood in which $f$ is one-to-one). Suppose that there is a positive
lower semicontinuous (l.s.c.) function $m:[0,\infty )\rightarrow \lbrack
0,\infty )$ such that condition (3.1) and the following condition are
satisfied%
\begin{equation*}
\underset{0}{\overset{\infty }{\dint }}m\left( s\right) ds=+\infty
\end{equation*}%
Then $f$ \ is a global homeomorphism, the inverse mapping $f^{-1}$ is
locally Lipschitz, and for every $y\in F$, the Lipschitz constant of $f^{-1}$
at $y$ \ is not greater than $m\left( \left\Vert f^{-1}\left( y\right)
\right\Vert \right) ^{-1}.$
\end{theorem}

In \cite{katriel} Katriel proved that from Ioffe's global inversion theorem
(that is from Theorem \ref{TheoFromIoffe} stated above) one can obtain the
Hadamard-Levy theorem. The proof of Theorem \ref{TheoFromIoffe} was based on
the Liusternik theorem:

\begin{theorem}
If $E,F$ are two Banach spaces and $f:E\rightarrow F$ is a local
diffeomorphism at $x$ then%
\begin{equation*}
\text{sur}\left( f,x\right) =\frac{1}{\left\Vert \left[ f^{\ \prime }\left(
x\right) \right] ^{-1}\right\Vert }
\end{equation*}
\end{theorem}

Of course one can easily see that from Ioffe's global inversion theorem one
can obtain Pourciau's global inversion theorem.

One interesting problem is the following: What is the relation between the
Hadamard-Palais theorem and the Hadamard-Levy theorem ? The answer is that
in the finite dimensional case the Hadamard-Palais theorem gives necessary
and sufficient conditions for global invertibility of local diffeomorphisms
but the Hadamard-Levy theorem gives only sufficient conditions. However
checking the coercivity condition from the Hadamard-Palais theorem is more
difficult than checking the integral divergence condition from the
Hadamard-Levy theorem.\bigskip

In the following subsection we shall give an example of a nonsmooth function
which verifies conditions from Hadamard-Palais theorem but do not satisfy
conditions from Pourciau theorem.

\subsection{\textbf{Example}}

For every $a\in \left( -1,1\right) $ consider the function $f_{a}:%
\mathbb{R}
^{2}\rightarrow 
\mathbb{R}
^{2}$ defined by

\begin{equation*}
f_{a}\left( x,y\right) =(x+a\left\vert x\right\vert ,x^{3}+y)
\end{equation*}

Then the following assertions hold:

1. For every $a\in \left( -1,0\right) \cup \left( 0,1\right) $ the function $%
f_{a}$ is locally Lipschitz and non-differentiable, while for $a=0$\ it is a 
$C^{1}$\ function

2. For every $a\in \left( -1,1\right) $ and\ any $\left( x,y\right) \in 
\mathbb{R}
^{2}$\ every (generalized Jacobian) matrix $A\in \partial f_{a}\left(
x,y\right) $ is nonsingular.

3. $\left\Vert f_{a}\left( x,y\right) \right\Vert \rightarrow \infty $ as $%
\left\Vert \left( x,y\right) \right\Vert \rightarrow \infty $

4. If $A$ is a matrix we denote $\left[ A\right] =\underset{\left\Vert
u\right\Vert =1}{\inf }\left\Vert Au\right\Vert $, $m\left( t\right) =%
\underset{\left\Vert z\right\Vert \leq \ t}{\inf }\left[ \partial
f_{a}\left( z\right) \right] =\underset{\left\Vert z\right\Vert \leq \ t}{%
\inf }\underset{\ \ \ A\in \partial f_{a}\left( z\right) }{\inf }\left[ A%
\right] $ \ \ then $\underset{0}{\overset{\infty }{\dint }}m\left( t\right)
dt<+\infty $

As a consequence we have the following

\textbf{Conclusion. }\textit{For every }$a\in \left( -1,1\right) $\textit{\
the function }$f_{a}$\textit{\ is a global diffeomorphism that satisfies
conditions from the Hadamard-Palais theorem, Theorem \ref{MainTheo}, and do
not satisfy the conditions from the Pourciau's global inversion theorem.}%
\bigskip

\textbf{Proof of the conclusion. }\bigskip

A.\textit{\ }$f_{a}$\textit{\ is locally Lipschitz with }$\partial
f_{a}\left( x,y\right) $\textit{\ containing only nonsingular matrices for
any }$\left( x,y\right) \in 
\mathbb{R}
^{2}.$\bigskip

Note that $f_{a}\left( x,y\right) =(x+a\left\vert x\right\vert
,x^{3}+y)=(x,x^{3}+y)+a(\left\vert x\right\vert ,0)=f_{1}\left( x,y\right)
+af_{2}\left( x,y\right) $.

Note that $f_{1}$ is a $C^{1}$ function, while $f_{2}$ is locally Lipschitz.
Note that $A\in \partial f_{a}\left( x,y\right) $ if $A=\left[ 
\begin{array}{ll}
1 & 0 \\ 
3x^{2} & 1%
\end{array}%
\right] +a\left[ 
\begin{array}{ll}
t & 0 \\ 
0 & 0%
\end{array}%
\right] $, where $t\in \left[ -1,1\right] $. Since $a\cdot t\in \left(
-1,1\right) $, we see that $\partial f_{a}\left( x,y\right) $ contains only
nonsingular matrices for any $\left( x,y\right) \in 
\mathbb{R}
^{2}$, and any fixed $a\in \left( -1,1\right) $.

\bigskip

B.\textit{\ }$f_{a}$\textit{\ is coercive.}\bigskip

In order to prove that the function $f_{a}$ is coercive it suffices to prove
that $f_{a}$ is coercive for the $l^{1}$ norm on $%
\mathbb{R}
^{2}.$ For every $a\in \left( -1,1\right) $ denote

\begin{equation*}
u_{a}\left( x,y\right) =\left\Vert f_{a}\left( x,y\right) \right\Vert
_{1}=\left\vert x+a\left\vert x\right\vert \right\vert +\left\vert
x^{3}+y)\right\vert ,\ \left( x,y\right) \in 
\mathbb{R}
^{2}.
\end{equation*}

Note that

\begin{equation}
u_{a}\left( x,y\right) \geq \left\vert x\right\vert -\left\vert a\right\vert
\cdot \left\vert x\right\vert +\left\vert x^{3}+y)\right\vert \geq \left(
1-\left\vert a\right\vert \right) \left( \left\vert x\right\vert +\left\vert
x^{3}+y)\right\vert \right) =\left( 1-\left\vert a\right\vert \right)
u_{0}\left( x,y\right) \   \label{3.2}
\end{equation}%
\ \ \ \ \ \ \ \ \ \ \ \ \ \ \ \ \ \ \ \ \ \ \ \ \ \ \ \ \ \ \ \ \ \ \ \ \ \
\ \ \ \ \ \ \ \ \ 

Consider the sets: 
\begin{equation*}
B_{1}=\left\{ \left( x,y\right) \in 
\mathbb{R}
^{2}:x\geq 0,y\geq 0\right\} ,
\end{equation*}

\begin{equation*}
B_{2}=\left\{ \left( x,y\right) \in 
\mathbb{R}
^{2}:x\leq 0,y\leq 0\right\} ,
\end{equation*}%
\begin{equation*}
B_{3}=\left\{ \left( x,y\right) \in 
\mathbb{R}
^{2}:x\geq 0,y\leq 0\right\} ,
\end{equation*}

\begin{equation*}
B_{4}=\left\{ \left( x,y\right) \in 
\mathbb{R}
^{2}:x\leq 0,y\geq 0\right\} ,
\end{equation*}%
\begin{equation*}
B_{31}=\left\{ \left( x,y\right) \in 
\mathbb{R}
^{2}:x\geq -\sqrt[3]{y}\geq 0,x\geq 1\right\} ,
\end{equation*}

\begin{equation*}
B_{32}=\left\{ \left( x,y\right) \in 
\mathbb{R}
^{2}:1\leq x\leq -\sqrt[3]{y}\right\} ,
\end{equation*}%
\begin{equation*}
B_{33}=\left\{ \left( x,y\right) \in 
\mathbb{R}
^{2}:0\leq x\leq 1\right\} .
\end{equation*}

Note that $B_{3}\subset B_{31}\cup B_{32}\cup B_{33}$ and $\overset{4}{%
\underset{i=1}{\cup }}B_{i}=%
\mathbb{R}
^{2}.$ We shall prove that for every $a\in \left( -1,1\right) $ the function 
$f_{a}$ is coercive on every $B_{i},i=1,2,3,4$ and on every $%
B_{3i},i=1,2,3.\bigskip $

If $\left( x,y\right) \in B_{1}$ then $u_{0}\left( x,y\right) =x+x^{3}+y\geq
x+y=\left\vert x\right\vert +\left\vert y\right\vert $, Note that

\begin{equation}
u_{a}\left( -x,-y\right) =u_{-a}\left( x,y\right)  \label{3.3}
\end{equation}%
hence $u_{a}\left( x,y\right) \geq \left( 1-\left\vert a\right\vert \right)
u_{0}\left( x,y\right) \geq \left( 1-\left\vert a\right\vert \right) \left(
\left\vert x\right\vert +\left\vert y\right\vert \right) $ \ for every $%
\left( x,y\right) \in B_{1}\cup B_{2}.\bigskip $

If $\left( x,y\right) \in B_{31}$ then 
\begin{equation*}
u_{0}\left( x,y\right) =x+x^{3}+y\geq x=\frac{1}{2}\left\vert x\right\vert +%
\frac{1}{2}\left\vert x\right\vert \geq \frac{1}{2}\left\vert x\right\vert +%
\frac{1}{2}\sqrt[3]{\left\vert y\right\vert }\geq \frac{1}{2}\sqrt[3]{%
\left\vert x\right\vert +\left\vert y\right\vert }
\end{equation*}

If $\left( x,y\right) \in B_{32}$ consider the function $\phi \left(
t\right) =t-t^{3}$, $t\in \lbrack 1,\infty ).$ Note that the function $\phi $
is decreasing and $\phi \left( x\right) \geq \phi \left( -\sqrt[3]{y}\right)
=-\sqrt[3]{y}+y.$

Then 
\begin{equation*}
u_{0}\left( x,y\right) =x-x^{3}-y=\phi \left( x\right) -y\geq -\sqrt[3]{y}%
+y-y=\frac{1}{2}\sqrt[3]{\left\vert y\right\vert }+\frac{1}{2}\sqrt[3]{%
\left\vert y\right\vert }\geq \frac{1}{2}\left\vert x\right\vert +\frac{1}{2}%
\sqrt[3]{\left\vert y\right\vert }\geq \frac{1}{2}\sqrt[3]{\left\vert
x\right\vert +\left\vert y\right\vert }
\end{equation*}

If $\left( x,y\right) \in B_{33}$ then 
\begin{equation*}
u_{0}\left( x,y\right) \geq \left\vert x\right\vert +\left\vert y\right\vert
-\left\vert x\right\vert ^{3}\geq \left\vert x\right\vert +\left\vert
y\right\vert -1
\end{equation*}

Since $u_{0}$ is coercive on $B_{31},B_{32}.B_{33}$ it follows that $f_{a}$
is coercive on $B_{3}.$ From (\ref{3.2}) and from the coercivity of $f_{a}$
on $B_{3}$ it follows the coercivity of $f_{a}$ on $B_{4}.\bigskip $

C. \textit{It holds that }$\underset{0}{\overset{\infty }{\dint }}m\left(
s\right) ds<+\infty $\textit{.}

For every $a\in \left( -1,1\right) $, $x\in 
\mathbb{R}
,t\in \left[ -1,1\right] $ we define the matrix

\begin{equation*}
Q\left( a,x,t\right) =\left( 
\begin{array}{cc}
1+at & 0 \\ 
3x^{2} & 1%
\end{array}%
\right)
\end{equation*}

Note that for every $a$ and $x\neq 0$ we have $\partial f_{a}\left(
x,y\right) =\left\{ Q\left( a,x,\frac{x}{\left\vert x\right\vert }\right)
\right\} .$

If $A=\left( a_{ij}\right) $ is a square matrix of order 2 then there exists
a constant $c_{1}>0$ such that $\left\Vert A\right\Vert \geq $ $c_{1}\left(
\left\vert a_{11}\right\vert +\left\vert a_{12}\right\vert +\left\vert
a_{21}\right\vert +\left\vert a_{22}\right\vert \right) $.

Note that if $x=0$ then

\begin{equation*}
\left[ Q\left( a,x,t\right) \right] ^{-1}=\left( 
\begin{array}{cc}
\left( 1+at\right) ^{-1} & 0 \\ 
0 & 1%
\end{array}%
\right)
\end{equation*}

\begin{equation*}
\left\Vert \left[ Q\left( a,x,t\right) \right] ^{-1}\right\Vert \geq
c_{1}\left( \left\vert \left( 1+at\right) ^{-1}\right\vert +1\right)
\end{equation*}

hence

\begin{equation*}
\underset{t\in \left[ -1,1\right] }{\sup }\left\Vert \left[ Q\left(
a,x,t\right) \right] ^{-1}\right\Vert \geq c_{1}\left( \left( 1+\left\vert
a\right\vert \right) ^{-1}+1\right)
\end{equation*}

Note that if $x\neq 0$ then

\begin{equation*}
\left[ Q\left( a,x,t\right) \right] ^{-1}=\frac{1}{1+a\frac{x}{\left\vert
x\right\vert }}\left( 
\begin{array}{cc}
1 & 0 \\ 
-3x^{2} & 1+a\frac{x}{\left\vert x\right\vert }%
\end{array}%
\right)
\end{equation*}

hence

\begin{equation*}
\left\Vert \left[ Q\left( a,x,t\right) \right] ^{-1}\right\Vert \geq c_{1}%
\frac{1+3x^{2}+\left\vert 1+a\frac{x}{\left\vert x\right\vert }\right\vert }{%
1+a\frac{x}{\left\vert x\right\vert }}\geq
\end{equation*}

\begin{equation*}
\geq c_{1}\left( 1+\frac{1+3x^{2}}{1+\left\vert a\right\vert }\right)
\end{equation*}

Note that

\begin{equation*}
\sup \left\{ \left\Vert \left[ Q\left( a,x,t\right) \right] ^{-1}\right\Vert
:\left\vert x\right\vert +\left\vert y\right\vert \leq s\right\} \geq
c_{1}\left( 1+\frac{1+3s^{2}}{1+\left\vert a\right\vert }\right)
\end{equation*}

hence

\begin{equation*}
m\left( s\right) =\underset{\left\Vert z\right\Vert \leq \ s}{\inf }\left[
\partial f_{a}\left( z\right) \right] \leq \frac{1+\left\vert a\right\vert }{%
c_{1}}\cdot \frac{1}{2+\left\vert a\right\vert +3s^{2}}
\end{equation*}

From the above inequality it follows that

\begin{equation*}
\underset{0}{\overset{\infty }{\dint }}m\left( s\right) ds<+\infty
\end{equation*}%
\hfill $\blacksquare $

\begin{tabular}{l}
Marek Galewski \\ 
Institute of Mathematics, \\ 
Lodz University of Technology, \\ 
Wolczanska 215, 90-924 Lodz, Poland, \\ 
marek.galewski@p.lodz.pl%
\end{tabular}%
\begin{tabular}{l}
Marius R\u{a}dulescu \\ 
Academia Romana \\ 
Institute of Mathematical Statistics \\ 
and Applied Mathematics \\ 
Calea 13 Septembrie nr 13, Bucharest, \\ 
RO-050711, Romania \\ 
mradulescu.csmro@yahoo.com;%
\end{tabular}

\end{document}